\date{March 11, 2019}

\documentclass[12pt]{amsart}
\usepackage{latexsym,amsmath,amsfonts,amscd,amssymb}
\usepackage{graphics}
\usepackage{xcolor}
\theoremstyle{remark}

\newcommand{\RR}{{\mathbb R}}

\newcommand{\ZZ}{{\mathbb Z}}

\title{Notes on the historical bibliography of the gamma function}


\author[R. P\'{e}rez-Marco]{Ricardo P\'{e}rez-Marco}
\address{CNRS, IMJ-PRG, Univ. Paris 7, Paris, France}
\email{ricardo@math.univ-paris13.fr}

\begin{document}

\maketitle

\begin{abstract}
Telegraphic notes on the historical bibliography of the Gamma function and Eulerian integrals. 
We correct some classical references. Some topics of the interest of the author. We provide 
some extensive (but not exhaustive) bibliography. Feedback is welcome, notes will be updated and 
some references need completion.
\end{abstract}

\section{Correction of references.}

\subsection{Integral formula for Euler-Mascheroni constant.}

Whittaker-Whatson (\cite{WW} 12.3 Example 2, p. 248) attributes the formula to Dirichlet (presumabibly to  \cite{Di}, 1836), 
but this formula  is already in Euler (1770) \cite{[E393]} section 25. Moreover, he also gives 
the formula in \cite{[E583]} (1785) and devotes a full article \cite{[E629]} (1789) to it.

\subsection{Frullani integral.}

According to Binet \cite{Bi} the integral is used before by Euler.

\subsection{Hadamard-Weierstrass product formula for the Gamma function.}

According to N\"orlund \cite{No5} and to Remmert p.41 \cite{Rem2},the product formula is due to Schl\"omilch \cite{Schl} p.171. 
According to Aycock \cite{Ayc3} p.6, the product formula was essentially found by Euler in Chapter 16 and 17 
of the second part of \cite{[E212]} and in \cite{[E613]}. Ph.J. Davis \cite{Da} attributes the product 
formula to F.W. Newman (1848) (no reference provided, but it is given in \cite{Go} and in \cite{Rem2} as \cite{Ne} p.57). 

Weierstrass theory of factorization is dated of 1876 \cite{Wei2}.

\subsection{Malmsten formula.}

Usually attributed to Malmsten (1849) \cite{Ma}, it is already in Binet (1839) \cite{Bi}, and 
N\"orlund \cite{No5} attributes the formula to Plana.

\subsection{Gauss multiplication formula.}

Apparently it is also already in Euler \cite{[E421]} according to A. Aycock (see \cite{Ayc2} for the details). 
Hermite knows that is due to Euler (see his cours \cite{He}, Le\c con 15, p.145). Artin attributes it to Gauss \cite{Ar} p.24.

\subsection{Laplace integral formula.}

Pribitkin article starts with the historically assertion that ``In 1812 Laplace establishes...'' 
citing \cite{Lap2} p.134 when indeed the formula was established 30 years earlier in 1782 in \cite{Lap1}.

\subsection{Kummer trigonometric expansion.}

Attributed to Kummer (1847)  but apparently found first by Malmst\'en (1846) (published in 1849 \cite{Mal},
but the submission date is 1st May 1846), according to \cite{Ayc1}.

\section{Different definitions of the Gamma function.}

\subsection{Bohr-Mollerup.}

Bohr and Mollerup \cite{BoMo} gave the characterization of the Gamma function as exercice in his calculus book. Emil Artin
made this the starting point of the definition of the Gamma function in his monograph \cite{Ar}. This 
definition became popular. Bourbaki,  \cite{Bou} Chapter VII, takes also 
this approach to define the Gamma function.
This definition is used and generalized by Vigneras \cite{Vi} to define higher Barnes Gamma functions.

\subsection{Wielandt's criterion.}

Based on the functional equation, plus estimate on vertical strip.
See Wielandt \cite{Wie} and Remmert \cite{Rem} (see Aycock \cite{Ayc3} p.12, see also Birkhoff \cite{Bir}).

\subsection{Weierstrass product.}

Attributed to Weierstrass in  Whittaker-Whatson \cite{WW}.

\subsection{Functional equation plus asymptotic estimates.}

See Prym \cite{Pr}. 

Weierstrass \cite{Wei} (see Aycock \cite{Ayc3} where he remarks that Weierstrass added the 
condition $\Gamma(1)=1$ which is not necessary in view of the estimate)  added the condition:
$$
\lim_{n\to +\infty} \frac{\Gamma(x+n)}{(n-1)!n^x}
$$

Usually the functional-difference equation is solved explicitly (see \cite{Ayc3} p.19, 
see also the Norlund approach). See also [Eu613] for the same ideas. Norlund fundamental solution 
for the difference equation \cite{No1} and \cite{No3} provides also a unique characterization 
of the Gamma function.

\subsection{Functional equation plus more symmetries.}

Davis \cite{Da} p.867 writes, without providing a reference,

\medskip

\textit{`` By the middle of the 19th century it was recognized that Euler's gamma function was the only 
continuous function which satisfied simultaneously the recurrence relationship, the reflection formula 
and the multiplication formula.''}

\subsection{Gauss product formula.} 

Product formula known to Euler. Gauss \cite{Ga}, Liouville \cite{Liou}.

\subsection{Integral representation.}

This is a very common approach in all modern analysis books and the original 
definition by Euler (1729-1730) who uses the integral:

$$
\Gamma(s) = \int_0^1 \left (\log (1/t)\right )^{s-1} dt
$$

The most popular version nowadays of integral form is:
$$
\Gamma(s)=\int_0^{+\infty} t^{s-1} e^{-t}\, dt
$$

Legendre, Liouville, Pringsheim \cite{Pr}.

\section{Various notes.}

\subsection{Weierstrass ``factorielle'' function.}

Apparently (\cite{Da} p.862), Weierstrass preferred to work with the function $Fc(u) =1/\Gamma(1+u)$ 
that he called ``factorielle'' function (probably in \cite{Wei}).

\subsection{Bourget function.}

See Godefroy \cite{Go} p.24. Bourget fonction $T(s)$ satisfies
$$
T(s+1)-T(s) =-\text{Fc}(s+1)
$$
where $\text{Fc}(s)$ is Weierstrass factorielle function.
Hence, $e^T(s)$ is $\Gamma_1^f(s)$ for $f(s)=e^{-\text{Fc} (s+1)}$.

We have 
$$
T(s)=\frac{eP(s)}{\Gamma(s)}
$$
where $P(s)$ is the Prym function
$$
eP(s)=\frac1s +\frac{1}{s(s+1)}+\frac{1}{s(s+1)(s+2)}+\ldots =
$$
also
$$
P(s)=\frac1s-\frac{1}{1}\frac{1}{s+1}+\frac{1}{1.2}\frac{1}{s+2}+\ldots +(-1)^n \frac{1}{n!}\frac{1}{s+n}+\ldots S
$$
The Prym function is the polar (Mittag-Leffler) part of the $\Gamma$ function:
$$
Q(s)=\Gamma(s)-P(s)
$$
is an entire function.

Problem mentioned in \cite{Go} p.23: Prove that $P$ has no more than 4 complex zeros. Not known if $Q$ has zeros.

\subsection{Stirling series.}

According to \cite{Da} p.862 (no reference provided), C. Hermite (1900) wrote down the Stirling series for $\log \Gamma(1+s)$,
convergent for $\Re s >0$,
$$
\log \Gamma(1+s) =\binom{s}{2}\log 2 +\binom{s}{3}(\log 3-2\log 2)+\ldots
$$

In the same place he also attributes to M.A. Stern (1847) (no reference provided) 
the Stirling series for the $\psi$ function, convergent for $\Re s >0$,

$$
\psi(s)= \frac{d}{ds}\log \Gamma(s) =-\gamma +\binom{s}{1}-\frac12 \binom{s}{2} 
+\frac13 \binom{s}{3}+\ldots
$$

\subsection{Hadamard interpollation of the factorial.}

According to \cite{Da} p.865, J. Hadamard (1894, no reference provided) 
gave an entire function interpollation of the factorial:

$$
H(s)=\frac{1}{\Gamma(1-s)} \frac{d}{ds} \left ( \frac{\Gamma\left 
(\frac{1-s}{2}\right )}{\Gamma \left (1-\frac{s}{2}\right )}  \right )
$$
It satisfies the functional equation
$$
H(s+1)=sH(s) +\frac{1}{\Gamma(1-s)}
$$

\subsection{Mellin Higher Gamma functions.}

According to Godefroy p.81, Mellin \cite{Me1897} defined general Gamma functions satisfying 
$$
F(s+1)=R(s)F(s)
$$
where $R$ is a rational function.

\subsection{Davis pseudo-Gamma function.}

Davis \cite{Da} p.867 gives the following pseudo-Gamma function $\Gamma_S$,
\begin{align*} 
\Gamma_S(s) &=1/s \ \ \ \ \text{for } 0<s<1 \\
\Gamma_S(s) &=1 \ \ \ \ \text{for } 1<s<2 \\
\Gamma_S(s) &=s-1 \ \ \ \ \text{for } 2<s<3 \\
\Gamma_S(s) &=(s-1)(s-2) \ \ \ \ \text{for } 3<s<4 \\
\ldots & \ldots \\
\Gamma_S(s) &= (s-1)^{\underline{k-1}} \ \ \ \ \text{for } k<s<k+1
\end{align*}

It is a convex function in $\RR_+^*$, and satisfies the functional equation
$$
\Gamma_S(s+1)=s\, \Gamma_S(s)
$$

\subsection{Bendersky Gamma function.}

Bendersky (1933, \cite{Be})  studies another hierarchy of Gamma functions different from Barnes'. Bendersky's Gamma functions 
have been rediscovered by Milnor 50 years later (1983, \cite{Milnor1983}) arising as higher partial derivatives of Hurwitz 
zeta function \`a la Lerch.


\begin{thebibliography}{9}

\bibitem{Ad} ADAMCHIK, V.; {\it Integrals associated with the multiple gamma function}, Integral Transforms and Special Functions, 
\textbf{25}, 6, p.462-469, 2014.

\bibitem{Ah} AHLFORS, L.V.; {\it Complex analysis}, 3rd edition, McGraw-Hill, 1979.

\bibitem{Al} ALEXEJEWSKY, W.; {\it Ueber eine Classe von Funktionen die der Gammafunktion analog sind}, Leipzig Weidmanncshe Buchhandluns, \textbf{46}, p.268-275, 1894.

\bibitem{Ayc1} AYCOCK, A.; {\it Translation of C.J. Malmst\`en's paper: ``De integrabilus quibusdam 
definitis, seriebusque infinitis''}, ArXiv:1306.4225, 2013.


\bibitem{Ayc2} AYCOCK, A.; {\it Euler and the multiplication formula for the $\Gamma$-function}, ArXiv:1901.03400, 2019.


\bibitem{Ayc3} AYCOCK, A.; {\it Euler and the Gammafunction}, ArXiv:1908.01571, 2019.


\bibitem{Ar} ARTIN, E.; {\it The gamma function}, Athena Series: Selected Topics in Mathematics, 
New York-Toronto-London: Holt, Rinehart and Winston, 1931.


\bibitem{Ba} BARNES, E.W.; {\it The theory of $G$-function}, Quat. J. Math., \textbf{31}, p. 264-314, 1900.



\bibitem{Ba2} BARNES, E.W.; {\it On the theory of the multiple Gamma function}, Cambr. Trans., \textbf{19}, p.374-425, 1904.



\bibitem{Be} BENDERSKY, B.C.; {\it Sur la fonction $\Gamma$}, Acta Mathematica, 1933, 263-322.



\bibitem{Bi} BINET, J.; {\it Sur les int\'egrales d\'efinies eul\'eriennes, et leur application \`a la 
th\'eorie des suites, ainsi qu'\`a l'\'evaluation de fonctions de grands nombres}, Journal 
de l'\'Ecole Polytechnique, \textbf{XVI}, p. 123-345, 1839.


\bibitem{Bir} BIRKHOFF, G.D.; {\it Note on the Gamma function}, Bull. Amer. Math. Soc., \textbf{20}, 1,  
p. 1-10, 1913.



\bibitem{BoMo} BOHR, H.; MOLLERUP, J.; {\it Laerebog i Kompleks Analyse}, vol. III, Copenhagen, 1922.



\bibitem{BC} BORWEIN, J.M.; CORLESS, R.M.; {\it Gamma and Factorial in the Monthly}, The American Math. 
Monthly, \textbf{121}, 1, 2017.


\bibitem{Bou} BOURBAKI, N.; {\it \'Elements de Math\'ematiques. Fonctions d'une variable r\'eelle}, Springer, 2006.



\bibitem{Da} DAVIS, Ph.J.; {\it Leonhard Euler's Integral: A Historical Profile of the Gamma Function}, 
The American Mathematical Monthly, \textbf{66}, 10, p. 849-869, 1959.



\bibitem{Di} DIRICHLET, P.G.L.; {\it Sur les int\'egrales eul\'eriennes}, Crelle, Journal f\"ur die reine und angewandte Mathematik, \textbf{15}, p. 258-263, 1836.

\bibitem{Eu} EULER, L.; {\it De Progressionibus harmonicus observatione}, 
Commentarii Academiae Scientarum Imperialis Petropolitanae, 7-1734, p.150-161, 1735.

\bibitem{[E00715]} EULER, L.; {\it Letter to Goldbach}, 13 October 1729,  
Euler Archive [E00715], eulerarchive.maa.org, 1729.

\bibitem{[E00717]} EULER, L.; {\it Letter to Goldbach}, 8 January 1730,  
Euler Archive [E00717], eulerarchive.maa.org, 1730.

\bibitem{[E212]} EULER, L.; {\it  Institutiones calculi differentialis cum eius usu in 
analysi finitorum ac doctrina serierum, volume 1}, Euler Archive [E212], eulerarchive.maa.org, 1787.


\bibitem{[E393]} EULER, L.; {\it De summis serierum numeros Bernouillianos involventium}, 
Commentarii Academiae Scientarum Imperialis Petropolitanae, 14, p.129-167; Opera Omnia, 
Series 1, \textbf{15}, p.91-130;  Euler Archive [E393], eulerarchive.maa.org, 1770.


\bibitem{[E421]} EULER, L.; {\it  Evolutio formulae integralis $\int x^{f-1} dx (lx)^{m/n}$ integratione a 
valora x=0 ad x=1 extensa}, Euler Archive [E421], eulerarchive.maa.org, 1771.


\bibitem{[E583]} EULER, L.; {\it  De numero memorabili in summatione progressionis harmonicae naturalis 
ocurrente}, Acta Academia Scientarum Imperialis Petropolitinae 5, \textbf{15}, p.569-603; 
Euler Archive [E583], eulerarchive.maa.org, 1785.


\bibitem{[E613]} EULER, L.; {\it  Dilucidationes in capita postrema calculi mei differentalis de
functionibus  inexplicabilibus}, Euler Archive [E613], eulerarchive.maa.org, 1787. Also, 
M\'emoires de l'acad\'emie des sciences de St. Petersbourg, \textbf{4}, 1813, p. 88-119; 
Opera Omnia: Series 1, \textbf{16}, p. 1-33.

\bibitem{[E629]} EULER, L.; {\it  Evolutio formulae integralis 
$\int dx \left ( \frac{1}{1-x} +\frac{1}{\log x}\right )$ a termino $x=0$ ad $x=1$ extensae}, 
Nova Acta Academiae Scientarum Imperialis Petropolitinae 4, p.3-16; Opera Omnia: Series 1, Volume 18, p.318-334; 
Euler Archive [E629], eulerarchive.maa.org, 1787.,


\bibitem{[E652]} EULER, L.; {\it De termino generali serierum hypergeometricarum}, 
Acta Academia Scientarum Imperialis Petropolitinae 7, p.42-82; 
Euler Archive [E652], eulerarchive.maa.org, 1793.

\bibitem{Fru} FRULLANI, G.; {\it Sopra gli integrali definiti}, Societa Ital. Scienze, \textbf{20}, p.44-48, 1828.


\bibitem{Fu} FUGLEDE, B.; {\it A sharpening of Wielandt's characterization of the Gamma function}, 
The American Math. Monthly, \textbf{115}, 9, p.845-850, 2008.


\bibitem{Ga} GAUSS, C.F.; {\it Disquisitiones generales circa series infinitam etc. Pars prior}, Werke, \textbf{III}, p.159, 1812.


\bibitem{Ga2} GAUSS, C.F.; {\it Letter to Bessel, December 18 1811}, Werke, \textbf{8}, p.90-92, 1880.



\bibitem{G2} GLAISHER, J. W. L.; \textit {On the history of Euler's constant}, Messenger Math., \textbf{1}, p.25-30, 1872.



\bibitem{Go} GODEFROY, M. ; {\it La fonction Gamma; Th\'eorie, histoire, bibliographie}, Gauthier-Villars, 1901.



\bibitem{Hankel} HANKEL, H. ; {\it Die Euler'schen Integrale bei unbeschr\"ankter Variabilit\"at des Argumentes}, Zeitschr. Math. Phys., 
\textbf{9}, p.1-21, 1864.


\bibitem{He} HERMITE, Ch. ; {\it Cours de M. Hermite}, r\'edid\'e par M. Andoyer, 1882, 4\`eme \'edition, Paris, 1891.


\bibitem{Je1916} JENSEN, J.L.W.V. (GRONWALL, T.H.); {\it An elementary exposition of the theory of the gamma function}, 
Annals of Mathematics, 2nd series, \textbf{17}, 3, p.124-166, 1916.


\bibitem{Ka} KAIRES, H.H. ; {\it On the optimality of a characterization theorem for the Gamma function 
using the multiplication formula}, Aequationes Mathematicae, \textbf{51}, p.115-128, 1996.


\bibitem{Ki} KINKELIN, H.; \textit{Ueber eine mit der Grammafunction verwandlte Transcendente und deren Anwendung auf die Integralrechnung}, Crelle, 
Journal f\"ur die reine und angewandte Mathematik, \textbf{57}, p.122-138, 1860.


\bibitem{La} LAGARIAS, J.C.; {\it Euler's constant: Euler's work and modern developments}, Bulletin 
of the American Math. Soc., \textbf{50}, 4,  p.527-628, 2013.


\bibitem{Lap1} LAPLACE, P.S.; {\it M\'emoire sur les approximations des formules qui sont des fonctions de tr\`es grands nombres}, 
Oeuvres, Vol. 10, p.209-291. M\'emoires de l'Acad\'emie des Sciences de Paris, 1782.

\bibitem{Lap2} LAPLACE, P.S.; {\it Th\'eorie analytique des probabilit\'es}, Paris, 1812.

\bibitem{Le} LEGENDRE, A.M.; {\it Exercises de calcul int\'egral}, Paris, 1811.


\bibitem{Liou} LIOUVILLE, ; {\it}, Comptes Rendus Acad\'emie des Sciences, XXXV, p. 320, 1853.

\bibitem{Mal} MALMST\`EN, C.J.; {\it De integralibus quibusdam definitis, seriebusque infinitis}, Crelle, Journal f\"ur die 
reine und angewandte Mathematik, \textbf{38}, p.1-39, 1849.


\bibitem{Ma} MASCHERONI, L.; {\it Adnotationes ad calculum integralem Euleri}, Vol. 1 and 2, Ticino, Italia, 1790, 1792.
(Reprinted in Euler, L. Leonhardi Euleri Opera Omnia, Ser. 1, Vol. 12. Leipzig, Germany: Teubner, pp. 415-542, 1915)


\bibitem{Me} MELLIN, HJ.; {\it Zur Theorie der Gammafunction}, Acta Mathematica, \textbf{8}, p. 37-80, 1886.

\bibitem{Me1897} MELLIN, HJ.; {\it \"Uber hypergeometrische Reihen h\"oherer Ordnungen}, Acta Societatis Scientiarum Fennicae, \textbf{23}, 7, p.
1897.


\bibitem{Milnor1983} MILNOR, J.; {\it On polylogarithms, Hurwitz zeta functions, and the Kubert identities}, 
L'enseignement Mathématique, \textbf{29}, p.281-322, 1983.



\bibitem{Ne} NEWMAN, F.W.; {\it On $\Gamma(a)$ especially when $a$ is negative}, 
The Cambridge and Dublin Mathematical Journal, \textbf{3}, 2, p.57-60, 1848. 

\bibitem{Nie} NIELSEN, N.; {\it Handbuch der Theorie der Gammafunktion}, Leipzig, reprinted 
Chelsea Publ. Co., New York, (1956), 1906. 



\bibitem{No1} N\"ORLUND, N.E.; {\it M\'emoire sur le calcul aux diff\'erences finies}, 
Acta Mathematica, \textbf{44}, Paris, 1922. 


\bibitem{No2} N\"ORLUND, N.E.; {\it M\'emoire sur les polynomes de Bernoulli}, 
Acta Mathematica, \textbf{43}, Paris, 1923. 


\bibitem{No3} N\"ORLUND, N.E.; {\it Vorlesungen \"uber Differenzen-Rechnung}, Verlag von Julius Springer, Berlin, 1924.


\bibitem{No4} N\"ORLUND, N.E.; {\it Le\c cons sur les s\'eries d'interpollation}, 
Gauthier-Villars, Paris, 1926. 


\bibitem{No5} N\"ORLUND, N.E.; {\it Le\c cons sur les \'equations lin\'eaires aux diff\'erences finies}, 
Gauthier-Villars, Paris, 1929. 


\bibitem{Pri} PRIBITKIN, W. de A.; {\it Laplace's integral, the Gamma function, and beyond}, The American Math. Monthly, 
\textbf{109}, 3, p.235-245, 2002.


\bibitem{Pr} PRINGSHEIM, A.; {\it Zur Theorie der Gamma-Functionen}, Math. Annalen, \textbf{31}, 1888.

\bibitem{Pr} PRYM, F.E.; {\it Zur Theorie der Gammafunction}, J. Reine Angew. Math, \textbf{82}, 1877



\bibitem{Rem} REMMERT, R. {\it Wielandt's theorem about the $\Gamma$-function}, 
The Amer. Math. Monthly, \textbf{103},3, 1996.

\bibitem{Rem2} REMMERT, R. {\it Classical topics in complex function theory}, Graduate Texts in Mathematics, \textbf{172}, Springer, 1998.


\bibitem{Sch} SCHL\"OMILCH, O.; {\it Einiges \"uber die Eulerische Integrale der zweiten Art }, 
Arch. Math. Phys., \textbf{4}, 4, p.167-174, 1843.


\bibitem{Schl} SCHL\"OMILCH, O. {\it Einiges \"uber die Eulerische Integraleder zweiten Art}, Arch. MAth. Phys., \textbf{4}, p.167-174, 1843.


\bibitem{Sm} SMITH, W.D.; {\it The Gamma function revisited}, 
Internet preprint, schule.bayernport.com/gamma/gamma05.pdf, 2006.


\bibitem{Sri} SRINIVASAN, G.K. {\it The Gamma function: An eclectic tour}, The American Math. Monthly, 
\textbf{114}, 4, p.297-315, 2007.


\bibitem{Sti} STIRLING, J.; {\it Methodus differentialis sive tractatus de summatione et interpola-
tione serierum infinitarum}, London, 1730.(see also the annotated translation by I. Tweedle, Springer, 2003). 


\bibitem{Vi} VIGN\'ERAS, M.F.; {\it L'\'equation fonctionelle de la fonction z\^eta de Selberg du groupe modulaire $\text{PSL}(2,\ZZ)$},  
Ast\'erisque, \textbf{61}, p.235-249, 1979.


\bibitem{Wei} WEIERSTRASS, K.; {\it \"Uber die Theorie der analytischen Fakult\"aten},  
Journal f\"ur Mathematik, \textbf{51}, p.1-60, 1856.

\bibitem{Wei2} WEIERSTRASS, K.; {\it Zur Theorie der endeutigen Analytischen Functionen},  Mathematische Werke, 
\textbf{II}, p.77-124, 1876. (See also, E. Picard translation, Ann. Sc. Ec. Norm. Sup., 2\`eme s\'erie, 
\textbf{8}, p.111-150, 1879).


\bibitem{Wei3} WEIERSTRASS, K.; {\it Vorlesungen \"uber die Theorie der elliptischen Funktionen},  Mathematische Werke, 
\textbf{5}, adapted by J. Knoblauch, 1876).

\bibitem{Wie} WIELANDT, H.; {\it },  Mathematische Werke, \textbf{2}, 2, De Gruyter, New York, published 1996. 
See also \cite{Rem}.


\bibitem{WW} WHITTAKER, E.T.; WATSON, G.N.; {\it A course in modern analysis}, Cambridge Univ. Press, 4th edition, 1927.

\end{thebibliography}
\end{document}